\documentclass[12pt]{amsart}
\usepackage{amscd}
\usepackage{verbatim}
\usepackage{amssymb, amsmath, amsthm, amscd,ifthen}
\usepackage[dvips]{graphics}
\usepackage[cp866]{inputenc}
\usepackage{graphicx}
\usepackage{epsfig}

\usepackage{amsmath,amssymb,amscd,}
\usepackage[all]{xy}

\usepackage[dvips]{graphics}
\usepackage[cp866]{inputenc}
\usepackage{graphicx}
\usepackage{epsfig}
\usepackage[hang]{subfigure}

\theoremstyle{plain}
\newtheorem{theorem}{Theorem}
\newtheorem{lemma}{Lemma}

\newtheorem{proposition}{Proposition}

\theoremstyle{definition}
\newtheorem{definition}{Definition}
\newtheorem{example}{Example}

\theoremstyle{plain}
\newtoks\thehProclaim
\newtheorem*{Proclaim}{\the\thehProclaim}

\theoremstyle{definition}
\newtoks{\thehRemark}
\newtheorem*{Remark}{\the\thehRemark}

\begin{document}

\title{On the area of constrained polygonal linkages}

\author{ Gaiane Panina, Dirk Siersma}

\address{G. Panina: St. Petersburg Department of Steklov Mathematical Institute,  gaiane-panina@rambler.ru;  D. Siersma: Utrecht University, Department of Mathematics, d.siersma@uu.nl}



\begin{abstract}
We study configuration spaces of  linkages whose underlying graph are  polygons with diagonal constrains, or more general, partial two-trees. We show that (with an appropriate definition) the
oriented area is a Bott-Morse function on the configuration space. Its critical points are  described and Bott-Morse indices are computed.

This  paper  is a generalization of analogous results for polygonal linkages (obtained earlier by G. Khimshiashvili, G. Panina, and A. Zhukova).
\end{abstract}

\maketitle

\section{Introduction}
A polygonal linkage
 is a  linkage whose underlying graph is a polygon (or, equivalently, a single-cycle graph).  One thinks of it as of a flexible polygon with rigid edges and revolving joints at the vertices whose ambient space is the Euclidean plane.
The idea of considering the oriented area   as a Morse function on its configuration space has already led to some non-trivial results: the critical points (or, equivalently, critical configurations) are easily describable, and there exists a short formula for the Morse index \cite{khipan}, \cite{khipan1}, \cite{panzh}, \cite{zhu}. In some further generalization \cite{PanPerfMorse} the oriented area proves to be an exact Morse function.

In the present paper we extend  the class of underlying graphs of linkages in such a way that it is possible
to introduce the oriented area with the same nice Morse-theoretical properties.

We start in Section \ref{section_crit_simplecase}  with  \textit{three-chain linkages.} It is our first example which is not a polygonal linkage.  By definition, a three-chain linkage is a patch of three chains, and therefore the underlying graph has three cycles. One of the cycles is distinguished: we consider its area $S$ as the function defined on the configuration space. Generically, it is a Bott-Morse function. We prove that the critical configurations are characterized by a combination of\textit{ cyclic} and \textit{aligned} conditions. We also give a formula for  Bott-Morse indices of critical points and critical components (Section \ref{SecMorseXXX}).

As an interesting illustration,  we describe a Hessian bifurcation and show that it amounts to a pitchfork bifurcation (Section \ref{SecSmallEx}).

Finally, we present the most general result which includes  polygonal linkages and three-chains (Section \ref{SecPTT}). Namely, we show  that critical points characterization and Bott-Morse index formula  extend to  \textit{partial two-trees.}
This class of graphs is well-studied and widely used in computer science   since many  algorithmic problems  may be solved much more efficiently for partial k-trees, and in particular for partial two-trees, see \cite{site}. Partial two-trees have a number of graph-theoretical characterizations; we make use of some of them. Also linkages with underlying partial two-tree graph have some specific properties, see \cite{GaoSitharam}; the present paper describes yet another one.

Further generalization of our theory (that is, beyond partial two-trees) needs more complicated criteria of critical points; we present an  example.

\section{Preliminaries and notation}\label{section_preliminaries}

A \textit{ linkage} is a graph together with length assignment to its edges.
  Its \textit{configuration space } is the space of all  planar realizations  of the linkage.  Its \textit{reduced configuration space } is the quotient of the configuration space  by  orientation-preserving isometries of the plane, that is, by translations and rotations.

\medskip

\textbf{A realizability convention } We   assume  throughout the paper that the edge lengths of a linkage are always realizable, so that the configuration space is non-empty.

  A configuration of a linkage is called \textit{aligned}  if it fits in a straight line.

\medskip

  Assume that $X\neq Y$  are two vertices of  a configuration $P$ of some  linkage.
  If  (1)  $P$ is a smooth point of the configuration space, and (2)  the gradient of the distance function $|XY|$ does not vanish at $P$, we say
  that \textit{ for the configuration $P$, the vertex  $Y$ moves freely with respect to }$X$.

   This can be reformulated in two ways:
    (1) the map from the configuration space to the plane which sends a configuration to the vector $\overrightarrow{XY}$  is a submersion, and (2)
the vertex $Y$ admits a  non-zero infinitesimal motion in any direction in the plane (with $X$ fixed).

 Clearly, this notion is symmetric: for a given configuration, $X$ moves freely with respect to $Y$ iff $Y$ moves freely with respect to $X$. An example follows in Lemma \ref{LemmaOpenChain}.

\subsection{Open chains}
An  \textit{open chain}  is a  linkage whose underlying graph is a path graph.
In robotics it is also called \textit{a robot arm}. We assume that its vertices are numbered, so  an open chain has a natural orientation.
For a configuration $(p_1,...,p_{n+1})$ of an open chain consider the distance  $D=|p_1p_{n+1}|$ as
 a Morse function on the reduced configuration space  $(S^1)^{n-1}$. Its critical points are described in \cite{Farber}:

\begin{theorem}\label{ThmMorseAligned}\cite{MilKap}\label{ThmOpenChain}\begin{enumerate}
                 \item Critical configurations of the distance function $D$ are \textit{aligned configurations}.
                 \item The Morse  index of a critical configuration equals $f-1$, where $f$ is 
                the number of \textit{forward} edges,   i.e.  
                 the number of (directed) edges $\overrightarrow{p_ip}_{i+1}$ which have the same direction as the vector     $\overrightarrow{p_1p}_{n+1}$.
               \end{enumerate}
\end{theorem}

Let us reformulate the first statement of the theorem in a way to be used in the paper:

\begin{lemma}\label{LemmaOpenChain} For any configuration of an open chain, the terminal point \newline $T=p_{n+1}$ moves freely with respect to the initial point $I=p_1 $ unless the configuration  is aligned. \qed
\end{lemma}

\subsection{Oriented area for polygonal linkages}

A \textit{polygonal linkage} is a linkage with a single-cycle underlying graph. We assume that all the vertices are numbered, so the cycle carries a natural orientation.

Throughout the paper we assume that no configuration fits  a straight line. This is equivalent to  smoothness of reduced configuration space of the polygonal linkage, see \cite{Farber}.
Let us be  more precise. The parameter space $\mathbb{R}^n=(l_1,...,l_n)$ is divided into \textit{chambers }
by hyperplanes (called ''walls'') of type $\sum_1^n\pm l_i =0$  (for all possible combinations of $\pm$'s). If $(l_1,...,l_n)$ lies on no wall, the reduced configuration space is smooth. For two $(l_1,...,l_n)$ belonging to one and the same chamber the reduced configuration spaces are diffeomorphic.

\begin{definition} \label{Dfn_area} The \textit{oriented area} of a polygon  with the vertices \newline $p_i = (x_i,
y_i)$  is defined by
$$2S(P) = (x_1y_2 - x_2y_1) + \ldots + (x_ny_1 - x_1y_n).$$
\end{definition}

\begin{definition}
    A
polygon  $P$  is  \textit{cyclic} if all its vertices $p_i$
lie on a circle.

\end{definition}

 Cyclic polygons  arise in the framework of
 the paper as critical points of the oriented area. The following fact was observed first by Thomas Banchhoff (unpublished):

\begin{theorem}\label{Thm_crirical_are_cyclic}\cite{khipan}  Generically, $S$ is a Morse function.
 At smooth points of the reduced configuration space, a polygon $P$ is a critical point of the
oriented area  $S$  iff $P$ is a cyclic configuration.
        \qed
\end{theorem}

We recall a short formula for Morse index of a cyclic configuration from \cite{zhu}, \cite{khipan1}.
We fix the following notation for a cyclic configuration.

$\omega_P$ is the winding number of $P$ with respect to the center of the circumscribed circle
$O$.

$\alpha_i$  is the half of the angle between the vectors
$\overrightarrow{Op_i}$ and $\overrightarrow{Op_{i+1}}$. The angle
is defined to be positive, orientation is not involved.

$\varepsilon_i$ is the
orientation of the edge $p_ip_{i+1}$, that is,

 $\varepsilon_i=\left\{
                       \begin{array}{ll}
                         1, & \hbox{if the center $O$ lies to the left of } p_ip_{i+1};\\
                         -1, & \hbox{if the center $O$ lies to the right of } p_ip_{i+1}.
                       \end{array}
                     \right.$

$e(P)$ is the number of positive entries in $\varepsilon_1,...,\varepsilon_n$.

$\mu_P=\mu_P(S)$ is the Morse index of the function $S$ in the point
P. That is, $\mu_P(S)$ is the number of negative eigenvalues of the
Hessian matrix $Hess_P(S)$.

\begin{theorem}\label{Thm_Morse_closed_plane}\cite{panzh}, \cite{khipan1} Generically,
for a cyclic configuration $P$ of a polygonal linkage,
$$\mu_P(S)=\left\{
       \begin{array}{ll}
         e(P)-1-2\omega_P &\hbox{if }\   \sum_{i=1}^n \varepsilon_i \tan \alpha_i>0; \\
         e(P)-2-2\omega_P & \hbox{otherwise}.\qed
       \end{array}
     \right.$$
\end{theorem}

\section{ Oriented area and its critical points for three-chains}\label{section_crit_simplecase}

Take three oriented open chains $A,B,$ and $Z$,  and glue together their initial points $I$, and also glue together their terminal points $T$. The obtained linkage is a \textit{three-chain}.
If $A$,$B$, and $Z$  have $p$, $q$, and $r$ edges respectively, we say that the three-chain linkage has the type $[p,q;r]$.

The underlying graph has three simple cycles: $AB$, $AZ$, and $BZ$.

Let us first look at the reduced configuration space of a three-chain.
\begin{theorem} If none of the three cycles  $AB$, $AZ$, and $BZ$ of a three-chain linkage has an aligned configuration,  the reduced configuration space of the three-chain is a smooth manifold.
\end{theorem}
Proof. We mimic the proof of the analogous statement for polygonal linkages, see \cite{Farber}. Let us remove the edge number $r$ from the chain $Z$. We obtain ''a polygon $AB$ with a tail'' whose reduced configuration space $M$ is smooth.
Next we define the function $D$ on $M$ as the (squared) distance between the end of the tail and the point $T$. The reduced configuration space of the three-chain is the level set of the function $D$: it equals $D^{-1}(l)$, where $l$ is the length of the removed edge.
It is a smooth manifold unless $l$ is a critical value of the function $D$. Lemma \ref{LemmaOpenChain} implies that this happens only if one of the cycles gets aligned. \qed

\medskip


In this way we get a set of linear equations of type:
$$ \epsilon_1l_1 +\cdots +\epsilon_n \l_n=0 , \mbox{where} \; \epsilon
\in \{-1,0,1\} $$
The corresponding arrangement of hyperplanes (''\textit{walls}'') divides  the parameter space into convex \textit{chambers}. In each chamber the reduced configuration spaces are smooth and mutually diffeomorphic.

In the sequel we will concentrate on three-chains with smooth reduced configuration spaces only, but in case of singularities several statements are also valid in the smooth part.

\bigskip

\begin{figure}[h]
\centering \includegraphics[width=10 cm]{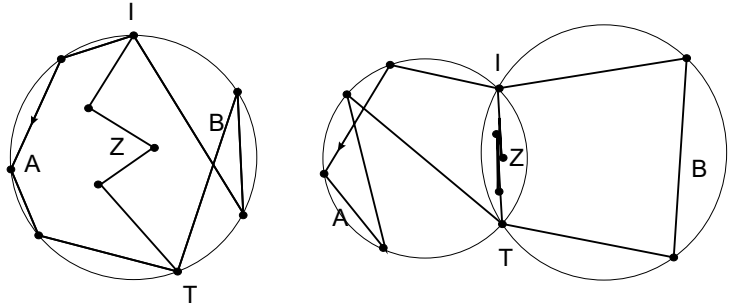}
\caption{Critical configurations of  a three-chain.}\label{crit}
\end{figure}

For a three-chain
we distinguish the oriented cycle $\Gamma=AB$  assuming that $A$ goes in the positive direction (and therefore, $B$ goes in the negative direction). In figures we indicate the orientation by an arrow, see Fig. \ref{crit}.

\begin{definition}
The \textit{oriented area of a three-chain linkage }is  the function  $S$ defined on the reduced configuration space as the
 the oriented area of the cycle  $\Gamma$.
 \end{definition}

\medskip

\textbf{A remark on genericity.} Throughout the paper we assume that the edge lengths are \textit{generic}.
We do not specify the meaning in advance, instead we keep a right to exclude from ''generic'' any finite number of manifolds of codimension at least one.
So the  generic case remains an open and everywhere dense subset of the parameter space.

\medskip

\begin{figure}
\centering
\includegraphics[width=12 cm]{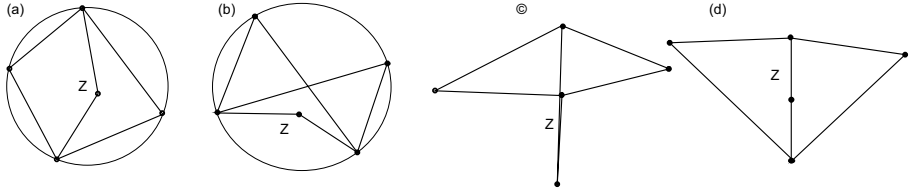}
\caption{Critical configurations: (a) and (b) are cyclic, (c) and (d) are aligned.  Each of them comes in a four-tuple.  }\label{Figure_4types}
\end{figure}

Let us first study the simplest case: three-chains of type $[2,2;2]$, or \textit{generalized Peaucellier-Lipkin linkage}. In this case the reduced configuration space is one dimensional, and therefore is a union of circles.
\begin{proposition}\label{LemmaSmall}  For a $[2,2;2]$ three-chain, critical points of $S$ fall into two types (see Figure \ref{Figure_4types}):

\begin{enumerate}
  \item Circular type: $\Gamma$ is realized as a cyclic polygon, that is, all the vertices of $\Gamma$ lie on a circle.
  \item Aligned type:  $Z$ is aligned.
\end{enumerate}
\end{proposition}

\begin{figure}[h]
\centering \includegraphics[width=5 cm]{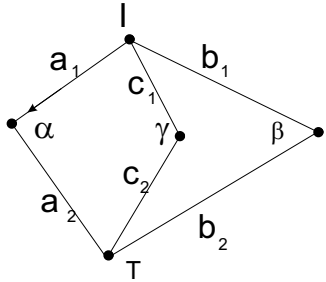}
\caption{Notation for Proposition \ref{LemmaSmall}.}\label{Notation1}
\end{figure}

Proof. In the notation of Figure \ref{Notation1} we have
$$2S=   a_1a_2 \cdot \sin \alpha+b_1b_2 \cdot \sin \beta .$$
We use the Lagrange multipliers   method with respect to the equations
$$a_1^2+a_2^2-2a_1a_2\cdot \cos \alpha=b_1^2+b_2^2-2b_1b_2\cdot \cos \beta$$

and
$$a_1^2+a_2^2-2a_1a_2\cdot \cos \alpha=c_1^2+c_2^2-2c_1c_2\cdot \cos \gamma.$$
The Lagrange multipliers  matrix is

$$M=\left(
    \begin{array}{ccc}
		a_1a_2\cdot \cos \alpha & b_1b_2\cdot \cos \beta & 0 \\
      2a_1a_2\cdot \sin \alpha & -2b_1b_2\cdot \sin \beta & 0 \\
      2a_1a_2\cdot \sin \alpha & 0 & -2c_1c_2\cdot \sin \gamma \\
    \end{array}
  \right)
.$$

At critical points its determinant vanishes. The direct computation gives
$$det M=  4 c_1c_2a_1a_2b_1b_2\cdot \sin \gamma \cdot \sin(\alpha+\beta),$$

which vanishes either $\sin\gamma =0$  (aligned type) or $\sin (\alpha +\beta)=0$  (circular type).\qed

\begin{theorem}  For a  three-chain linkage all critical points of the function $S$ fall into two  (generically disjoint) types, see Figure \ref{crit}:
\begin{enumerate}
 \item Circular type: $\Gamma$ is realized as a cyclic polygon, that is, all the vertices of $\Gamma$ lie on a circle.
 For each cyclic configuration of $\Gamma$ we have a critical manifold which equals the reduced configuration space of a flexible polygon  obtained from $Z$ by attaching the corresponding diagonal of the cyclic polygon.
  \item Aligned type: $Z$ is aligned, and each of $A$ and $B$ are cyclic. The two circumscribed circles generically do not coincide.
      \end{enumerate}
\end{theorem}
The proof
 comes from using Proposition \ref{LemmaSmall} sufficiently many times. The ruling idea is to freeze some of the joints of the linkage and to get a small linkage treated in Proposition \ref{LemmaSmall}.  As an example, if $A$ and $B$ contain two bars each, and
$Z$ contains three bars, we freeze one of the joints of $Z$ and use Proposition \ref{LemmaSmall} for the first time, then freeze the
other joint of $Z$ and use Proposition \ref{LemmaSmall} for the second time.\qed

\bigskip

\textbf{A genericity remark.} According to our genericity convention, we may assume that
cyclic  polygon $AB$ and simultaneously aligned $Z$ never occur.
We may also assume that each of the cyclic polygons that appear in the above theorem is a non-degenerate Morse point (considered as a single polygon).

 \example{Three-chain [2,2;2].}

The reduced configuration space in the smooth case  consists of either two or four circles, depending on intervals between minimum and maximal distances between the endpoints of the arms $A,B,$ and $Z$.
  The global picture is very sensitive for the mutual position of these intervals  and gives rise to many different cases.
There exist $[2,2;2]$ three-chains with 16 critical points (which is the maximum possible number of critical points).
On two components of its reduced configuration space we have only one maximum and one minimum of aligned type. The other components contain four circular and two aligned critical points.

\example{Three-chains  [2,2;2] are [2,2;3].}
These are the first cases after [2,2;2]. Their  reduced configuration spaces are one and the same  two dimensional manifold. One can meet here the first Bott-Morse singularities:\\
The three-chain [2,2;3] has isolated aligned singularities and Bott-Morse if the 4-cycle is circular.\\
The three-chain [2,3;2] has only isolated singularities; both circular (5-gon) and aligned.\\


\section{Bott-Morse indices for the three-chain linkage}\label{SecMorseXXX}

We start with a reminder about non-isolated critical points, see \cite{BM}. A function $f:M^n\rightarrow \mathbb{R}$ is called \textit{Bott-Morse }  if for every component $C$ of its critical set $f$ is locally of the form $f=f_o-(x_1^2+...+x_{\lambda}^2)+(x_{\lambda+1}^2+...+x_k^2)$.  $\lambda=\lambda (C)$ is called the Bott-Morse index of the component.

As a consequence, a critical component is a smooth manifold of codimension $n-k$.

Standard formula from Bott-Morse theory is
$$\sum(-1)^{\lambda}\chi(C)=\chi(M^n).$$

\bigskip

Consider first the critical point of cyclic type. For a fixed cyclic configuration of $\Gamma$  critical points form a manifold which  is assumed to be already known: it is the reduced configuration space of a polygon formed by $Z$ with an extra bar $W=IT$.

\begin{theorem}
The Bott-Morse index of the critical manifold of the circular type equals the Morse index of the cyclic polygon, see Theorem \ref{Thm_Morse_closed_plane}.
\end{theorem}
Proof. The area function can be described as $S=S(\xi,\eta)$, where $\xi$ are coordinates for $\Gamma$,
and $\eta$  are coordinates for $Z$. Clearly, $S$ does not depend on $\eta$.\qed

\medskip

Now assume that we have a critical configuration of aligned type. The cycle $\Gamma$ splits in the homological sum of two cycles $AW$ and $BW$  obtained by by cutting $\Gamma$  along the diagonal $W= IT$. These two polygons are cyclic, and therefore can be treated as critical points of the area function on the reduced configuration space of polygonal linkages defined by their edge lengths.

Before we formulate a theorem let us fix some notation.
We introduce the vector $\overrightarrow{W}= \overrightarrow{IT}$.

\begin{theorem}\label{ThmMorseThreeCh} The Morse index of an aligned critical configuration is

$$\mu =\mu_A +\mu_B + \nu,$$

where
$\mu_{A} $  and  $\mu_{B} $ are the Morse indices of the polygons $AW$ and $BW$ obtained by adding the diagonal $W$ to $A$ and $B$ respectively.
The formula for $\mu_{A} $  and  $\mu_{B} $ is contained in Theorem \ref{Thm_Morse_closed_plane}.

$$\nu= \left\{
          \begin{array}{ll}
            $f-1$, & \hbox{if the pair of vectors $(\overrightarrow{W},\overrightarrow{O_AO_B})$ is positively oriented;} \\
            $r-f$ , & \hbox{otherwise,}
          \end{array}
        \right.$$
where $r$ the number of edges in the chain $Z$, $f$   is the number of forward edges in $Z$ and $O_A$, and $O_B$ are the centers of circumscribed circles.
\end{theorem}

Let us first
 prove the following lemma.
\begin{lemma}\label{lemmaSecondDer}
Assume that in a triangle we fix  side lengths $a$ and $b$ and let $c$ be variable.
 Let $S=S(c)$ be the area of the triangle.\footnote{We mean here the ''usual'' area, which is positive.}

Then $$\frac{\partial S(c)}{\partial c }=\left\{
                                                              \begin{array}{ll}
                                                                 |OM|, & \hbox{if $\gamma < \frac{\pi}{2}$;} \\
                                                                  -|OM|, & \hbox{otherwise.}
                                                                \end{array}
                                                              \right.
$$

Here $\gamma$ is the angle opposite to $c$,
$M$ is the center of the edge $c$, and $O$ is the center of the circumscribed circle (with radius $R$).
\end{lemma}
Proof of the lemma. First note, that in terms of $\gamma$ the area is strictly increasing between $0$ and $\frac{\pi}{2}$, maximal at  $\frac{\pi}{2}$ and decreasing between  $\frac{\pi}{2}$ and $\pi$.  According to the cosine rule:
$$c^2 = a^2+ b^2 - 2ab \cdot \cos \gamma \; \; \hbox{\rm so} \; \; 2c = 2ab \cdot \sin \gamma \cdot \frac{d \gamma}{dc}.$$

Moreover $2S = ab \cdot \sin \gamma$ ; therefore
$$\frac{\partial (2S)}{\partial c} = ab \cdot \cos \gamma \cdot \frac{d \gamma}{dc}= \frac{ c \cos \gamma}{\sin \gamma} = \frac{2R \cdot |OM|}{R} = \pm 2 |OM|. $$
\qed

\medskip

Now we are ready to prove the theorem.
In  the neighborhood of an  aligned critical point we choose the following coordinate system. It consists of three families of coordinates.  Namely, we
first choose coordinates $\alpha =(\alpha_1,...,\alpha_{p-2})$ for the polygon  $AW$ and
 $\beta=(\beta_1,...,\beta_{q-2})$ for the polygon $BW$ assuming that $w = |IT|$ is constant. We specify these coordinates in more detail, when necessary (they could be angles or lengths of diagonals, etc).
Next we choose coordinates $\gamma = (\gamma_1,..., \gamma_{r-1})$ for the chain  $Z$. Note that $\alpha,\beta,w$ are coordinates for the polygon $AB$.
Consider a part of the Taylor expansion of $S=S(\alpha,\beta,w)$:
$$S=S_0 + \Big(\frac{\partial S}{\partial w}\Big)(w-w_0) +
\frac{1}{2} \sum  \Big(\frac{\partial^2S}{\partial \alpha_i \partial \alpha_j}\Big)\alpha_i \alpha_j +
\frac{1}{2} \sum  \Big(\frac{\partial^2S}{\partial \beta_i \partial \beta_j}\Big)\beta_i \beta_j
$$
Note that $\frac{\partial^2S}{\partial \alpha_i \partial \beta_j} = 0.$
By Morse lemma and Theorem \ref{ThmOpenChain}, at the critical point we can suppose that:  \\
$w-w_0 =  -(\gamma_1^2 +\cdots \gamma_N^2 ) + (\gamma_{N+1}^2 +\cdots \gamma_{r-1}^2 ) $. After substituting in the formula this gives the 2-jet of S as function of $\alpha,\beta,\gamma$.
So the Hessian matrix has a block structure:

$$Hess (S)=\left(
             \begin{array}{ccc}
               H_A & 0 & 0 \\
               0 & H_B & 0 \\
               0 & 0 & \mathcal{A }\\
             \end{array}
           \right)
$$
The first block $H_A$ (respectively, $H_B$)  equals the Hessian matrix of the area of the polygon $AW$ (respectively, $BW$).
The third block is a diagonal matrix $\mathcal{A}$, which  up to the constant $(\frac{\partial S}{\partial w})$ equals the Hessian matrix of the distance function  $D$ of the open chain $Z$, see Theorem \ref{ThmOpenChain}.

To determine the sign of $(\frac{\partial S}{\partial w})$, we
 specialize the choice of coordinates $\alpha$ and $\beta$ as follows:
triangulate the polygon $AW$  (respectively, $BW$), and take the lengths
of the  diagonals $d^A_1,...,d^A_k$ (respectively, $d^B_1,...,d^B_m$) as coordinates.

Let   $S_1^A,...,S_{k-1}^A$ and $ S_1^B,...,S_{m-1}^B$  be the (oriented) areas of the triangles of the triangulations.
 We may assume that $S_1^A$  and  $S_1^B$  correspond to triangles that are incident to the edge  $W$.
Then $S=S_1^A+...+S_{k-1}^A+S_1^B+...+S_{m-1}^B$, however, only two of the summands depend on $w$. Eventually, Lemma \ref{lemmaSecondDer} implies:

$$(*) \ \ \ \Big(\frac{\partial S}{\partial w}\Big) = \Big(\frac{\partial S^A_1}{\partial w}\Big) + \Big(\frac{\partial S^B_1}{\partial w}\Big)  =
\left\{
\begin{array}{lll}
|O_AO_B|, & \hbox{if the pair of vectors $(\overrightarrow{W},\overrightarrow{O_AO_B})$ } \\
\  & \hbox{ is positively  oriented;} \\
- |O_AO_B|, & \hbox{otherwise.}
\end{array}
\right.
$$

For generic parameter values $S$ will be Morse and the  Morse index (which equals the number of negative eigenvalues) sums up.
\qed

\section{Hessian bifurcation} \label{SecSmallEx}

We discuss a bifurcations  occurring in generic one-parameter families at an aligned singularity.

\subsection*{Hessian for [2,2;2]}\label{Hes222}
We discover a pitchfork bifurcation when the Hessian vanishes. We assume an aligned critical point
and  use the notations from Figure \ref{Notation1}.  Let $w=|IT|$  be the length of the diagonal. In the aligned situation we can assume that $w= w_0 \pm \epsilon^2$ (where $\epsilon$ is a  coordinate for the chain $Z$ and the sign $\pm$ depends on the type of alignment.

Next consider the 2-jet of $S$ at $w_0$ in terms of $w$:

$$S=S_0 + \Big(\frac{dS}{dw}\Big)(w-w_0) + \frac{1}{2} \Big(\frac{d^2S}{dw^2}\Big)(w-w_0)^2 $$
So the 4-jet of $S$ at $w_0$ in terms of $\epsilon$ is :
$$S=S_0 + \pm \Big(\frac{dS}{dw}\Big) \epsilon^2 + \frac{1}{2} \Big(\frac{d^2S}{dw^2}\Big)\epsilon^4 $$

The Hessian is already computed in Lemma \ref{lemmaSecondDer} for the two triangles.
 (*) implies that
the Hessian is equal to zero if all the four points of the cycle are on one circle.
So we have now a combination of the circular and the aligned condition.
The formula for $S$ starts in this case with a term of degree four. This means that locally we have a singularity of type $A_3$. By changing the side length  slightly it is possible to unfold the singularity such that  the sign of the Hessian changes exactly at the initial value of the length vector.
The unfolding is of type $S=S_0 + \epsilon^4 \pm \lambda \epsilon^2$ and this describes  a pitchfork bifurcation: one aligned maximum splits into one aligned minimum and two maxima of circular type. By symmetry, we could also have  a similar bifurcation where maximum and minimum are interchanged.

\begin{figure}[h]
\centering \includegraphics[width=10 cm]{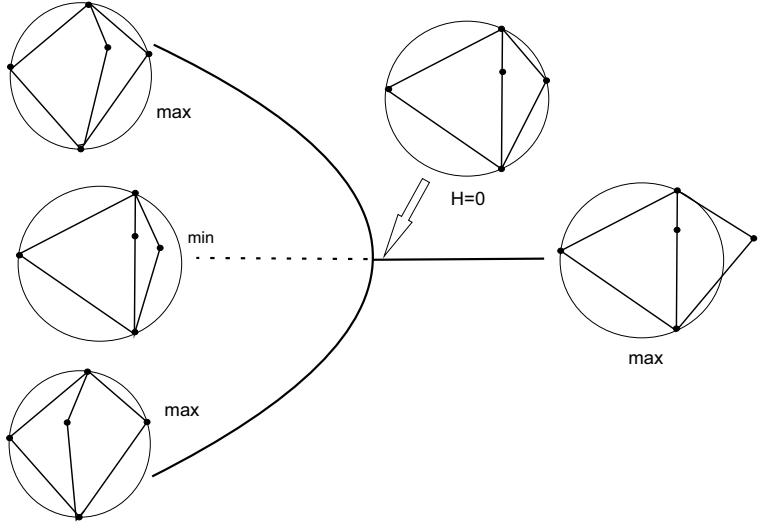}
\caption{Pitchfork bifurcation in the $[2,2;2]$ case.}\label{FigPitchfork}
\end{figure}

Note that this bifurcation takes place in the interior of a  chamber: during the deformation one still has diffeomorphic reduced configuration spaces.

This  pitchfork bifurcation also shows that in several cases our Morse function is not exact.

\subsection*{Hessian bifurcation for [p,q;r]}

Assume that $(\frac{\partial S}{\partial w})=0$. Then the Hessian becomes zero.
We take  the next term to the Taylor development of $S$  and get
$$S=S_0 + \Big(\frac{\partial S}{\partial w}\Big)(w-w_0) +\frac{1}{2} \Big(\frac{\partial^2 S}{\partial w^2}\Big)(w-w_0)^2 + Q(\alpha,\beta),
$$
where $Q(\alpha,\beta)$ is supposed to be a non-degenerate quadratic form. Remind that
$w-w_0 =  -(\gamma_1^2 +\cdots \gamma_s^2 ) + (\gamma_{s+1}^2 +\cdots \gamma_{r-1}^2 ), $.
so we work  essentially with degree 4 in $\gamma$-coordinates.

Note that passing through Hessian zero changes the sign in the matrix $\mathcal{A}$ and therefore the change in index can be greater than 1.

Let us consider the [2,2;3] case:
if $w-w_0 =  -\gamma_1^2 +\gamma_2^2 $ then there
  occur non-isolated singularities; the singular set $-\gamma_1^2 +\gamma_2^2 = 0$ is even not smooth.
Especially interesting is the case where $w-w_0 =  -(\gamma_1^2 +\gamma_2^2 )$. If the Hessian is zero there is an aligned isolated maximum  and this bifurcates from an isolated maximum into an isolated minimum together with a maximum of circle type.  This is a generalized pitchfork bifurcation.

\section{Area function for partial two-trees  and beyond}\label{SecPTT}

Let us start this section with an introductory example which lies in between three-chains and partial two-trees.

\begin{definition}
 A \textit{polygon with with non-crossing diagonals }is a
  graph $G$ with a distinguished cycle $\Gamma$  such that
 \begin{enumerate}
   \item $G$ can be embedded in the plane in such a way that
 $\Gamma$ is the outer cycle, and
   \item $G$ is obtained by attaching to $\Gamma$ a number of path graphs. Each path graph is attached by its initial and terminal points $I$ and $T$.
 \end{enumerate}
 In other words, $G$ is polygon with a number of non-crossing diagonals added. Each diagonal is a path graph.
\end{definition}
The above definition is combinatorial: we have defined a class of graphs, no edge lengths are involved so far.
So the words ''non-crossing diagonals'' have the combinatorial meaning only. For a configuration of an associated linkage
diagonals might cross.

\medskip

\begin{figure}[h]
\centering
\includegraphics[width=10 cm]{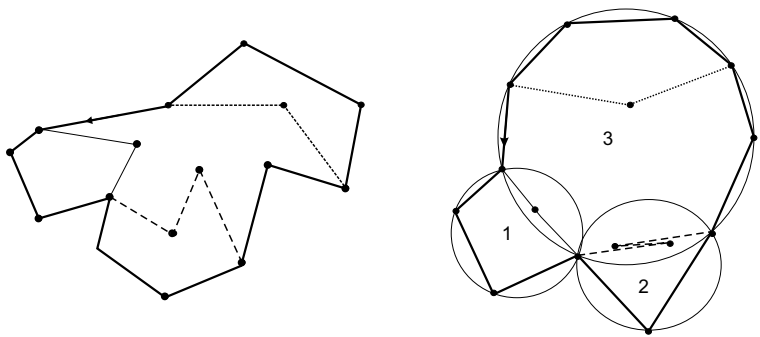}
\caption{A polygon with non-crossing diagonals and one of its critical configurations. Its orientation is indicated by an arrow.}\label{Figure_big}
\end{figure}

Next, we consider  linkages based on this type of graphs.
 It is easy to see that the approach of Sections \ref{section_crit_simplecase} and \ref{SecMorseXXX}  can be repeated literally (we give the statements in the full generality below in the section). Figure \ref{Figure_big}  depicts a polygon with non-crossing diagonals and one of the critical configurations.

A natural question is: how far can we extend the class of graphs keeping  a similar  circular-aligned characterization of critical points and easily computable Bott-Morse indices.
To our opinion, a proper generalization  are linkages whose underlying graphs are partial two-trees (PTTs).
It has been observed during AIM meeting ''Configuration spaces of linkages''
(October  2014)
 that for a PTT linkage, the parameter space  (that is, the space of the edge lengths of the linkage)
is divided by some hyperplanes (called '' walls'') into chambers. The diffeomorphic type of the reduced configuration space depends on the chamber only.  To the best of our knowledge, the result is not published yet. However, we do not use it here: we only  need that generically, the reduced configuration space of a partial two-tree is a smooth manifold.

\bigskip

A \textit{two-tree} is a graph obtained by starting with $K_3$ (complete graphs on three vertices) and then repeatedly taking  the connected sum  via edges with  a number of copies of $K_3$.

A \textit{partial two-tree} (PTT, for short) is a graph that can be characterized by the following equivalent conditions:
\begin{enumerate}
\item A partial two-tree becomes a two-tree after adding some extra edges.
  \item A partial two-tree has no $K_4$ minor, see \cite{WaldCoul}.

  \item
   PTTs are \textit{two-terminal series-parallel} graphs, see \cite{Epp}.

  Let us explain this in details.
  A \textit{two-terminal graph} (TTG) is a graph with two distinguished vertices, $I$ and $T$.
The \textit{parallel composition } of two TTGs $X$ and $Y$ is a a graph obtained by patching $I_X$ to $I_Y$, and also patching $T_X$ to $T_Y$.
The \textit{series composition}   of two TTGs $X$ and $Y$ is a a graph obtained by patching $I_Y$ to $T_X$.
Finally, a \textit{two-terminal series-parallel graph} is a two-terminal graph that may be constructed by a sequence of series and parallel compositions starting from a set of copies of a single-edge graph $K_2$.

\end{enumerate}

\bigskip

Given a PTT one can turn it to a two-terminal series-parallel by an appropriate choice of $I$ and $T$.
Not any two distinct vertices can serve as $I$ and $T$, so we shall need the following sufficiency condition:\footnote{Probably this lemma has already appeared in the literature.}

\begin{lemma}\label{LemmaComb1} If two vertices of a PTT are connected by an edge, they can serve as $I$ and $T$.
\end{lemma}
The proof goes by induction. We can assume that our PTT is a two-tree.  For two  vertices $Q$ and $R$ sharing an edge, two cases are possible: (1) $QR$ is incident to more that one triangle of the two-tree. Then $G$ is a connected sum of several smaller PTT over the edge $QR$, so we use the inductive assumption. (2) $QR$ is incident to exactly one triangle. Then removal of the edge $QR$ leaves a  a connected sum of two smaller PTT  over a single vertex incident to both $R$ and $Q$, and again inductive assumption is applicable.
\qed

\medskip

Assume now that we have a  linkage whose underlying graph is a partial two-tree.
An analog of Lemma \ref{LemmaOpenChain} is valid:

\begin{lemma}\label{LemmaPTT} For a configuration of a PTT linkage,  the vertex $T$ moves freely with respect to the vertex $I$  unless there exists an aligned path connecting $T$ and $I$.
\end{lemma}
The proof comes by induction based on the  above  characterization (3) combined with Lemma \ref{LemmaOpenChain}.
\qed

\medskip

Given a PTT linkage with underlying graph $G$, choose a cycle $\Gamma$ with no repeating vertices and take its area $S$ as the
function on the reduced configuration space of the linkage.

\medskip

We need the following combinatorial lemma.

\begin{lemma} \label{LemmaComb2}In the above notation  $G$ equals the connected sum of $\Gamma$ and a number of  disjoint graphs $G_i$.
\begin{enumerate}
\item Each of $G_i$ is a PTT, and therefore, a two-terminal series-parallel graph.
  \item The connected sum with each of $G_i$ is taken at at most two vertices, and
 \item  these two vertices $I_i$ and $T_i$ can serve as $I$ and $T$ for $G_i$.

\end{enumerate}
\end{lemma}
Proof. (1) is clear. (2) is true since $G$ has no $K_4$ minors. (3) follows from Lemma \ref{LemmaComb1}.
\qed

\medskip

We give the characterization of critical points of $S$ and the formula for the Bott-Morse indices in an algorithmic way:

\begin{theorem}\label{ThmPTTCrit}\textbf{Critical points and Bott-Morse indices for partial two-trees.}
Given a configuration $P$ of a PTT  linkage,
\begin{enumerate}
  \item take all aligned paths with endpoints on $\Gamma$ and replace them by single segments  which we call  \textit{straight line diagonals}. Since there is no $K_4$ minor, they never cross (in the combinatorial sense).
  \item Take the new configuration composed of  $\Gamma$ and straight line diagonals (the underlying graph changes!).
  \item The straight line diagonals decompose $\Gamma$ into the homological sum of elementary cycles $$\Gamma=\Gamma_1+...+\Gamma_k.$$
\end{enumerate}
Then we have:
\begin{enumerate}
  \item $P$ is a critical point of the area $S$ if and only if
each of $\Gamma_k$ has a circumscribed circle.

  \item For the Bott-Morse index of a critical point (or critical manifold) we have:

  $$\mu =\sum_i \mu_i + \sum_j\nu_j,$$
  where $\mu_i$ is the Morse index of $\Gamma_i$;
 $ \nu_j$refer to the $j$-th  aligned path.
  \item  All critical manifolds are products of reduced configuration spaces of partial two-trees with a smaller number of edges.
\end{enumerate}

\end{theorem}
Proof. (1) By Lemma \ref{LemmaComb2},  the vertices of $\Gamma$ move freely with respect to each other unless they are connected by an edge of $\Gamma$ or by an aligned path. Since on this step  it suffices to consider only first-order terms, we may think that  each aligned diagonal of $\Gamma$ fixes the distance between its endpoints.  The area of $\Gamma $ is critical if and only if all the $\Gamma_i$  are critical, and therefore all of them are cyclic.

 The proof of (2) follows the pattern of the  proof of Theorem \ref{ThmMorseThreeCh}.
Namely, after introducing appropriate coordinates the Hessian matrix becomes a block matrix. These are blocks that correspond to $\Gamma_i$ and to the aligned diagonals.  (3) is clear. \qed

\bigskip

\textbf{Example}

For the Morse index of the critical configuration depicted in Figure \ref{Figure_big}, we have:
$$\mu=\mu_1+\mu_2+\mu_3 +\nu_{1,3}+ \nu_{2,3} = $$$$1+0+5+(2-1)+(3-2)=8.$$

\medskip

\textbf{A non-example}

\begin{figure}[h]
\centering
\includegraphics[width=4 cm]{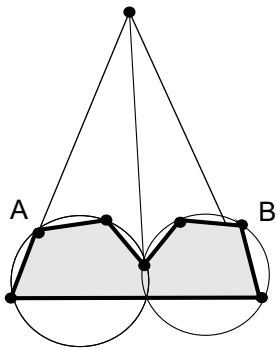}
\caption{A non-PTT. The cycle $\Gamma$ is designated bold.}\label{Figure_NonEx}
\end{figure}

The following example shows that Theorem \ref{ThmPTTCrit}  is not valid for arbitrary (non-PTT) graphs. Figure \ref{Figure_NonEx}
depicts an example. For this particular example it can be easily proven that critical configuration  must have cyclic $A$ and cyclic $B$. However, there is a one-parametric family of such configurations, with nonconstant  $S$ along the family. So being critical must have some other condition which is neither cyclic nor aligned.

\section*{Acknowledgments} This work is
 supported by the RFBR grant 17-01-00128.
It is our pleasure to acknowledge the  hospitality and excellent working conditions of
CIRM, Luminy, where
this paper was initiated as a 'research in pairs' project. We are grateful to George Khimshiashvili for useful remarks and comments.

\end{document}